\documentclass[12pt,a4paper]{article}
\usepackage[T2A]{fontenc}
\usepackage[utf8]{inputenc}
\usepackage[russian,english]{babel}
\usepackage{indentfirst}
\usepackage{misccorr}
\usepackage{graphicx}
\usepackage{amsmath}
\usepackage{amsthm}
\usepackage{thmtools}
\usepackage{longtable}
\usepackage{geometry}
\usepackage{imakeidx}
\usepackage{hyperref}
\usepackage{mathrsfs}
\usepackage{amssymb}
\usepackage{euscript}
\usepackage{csquotes}
\usepackage[backend=biber,style=numeric,sorting=none]{biblatex}
\newtheorem{definition}{Definition}
\newtheorem{theorem}{Theorem}

\DeclareMathOperator{\ran}{Ran}
\DeclareMathOperator{\dom}{Dom}

\newcommand{\W}[2]{\mathbb{W}_{\mathbb{#1}#2}}
\newcommand{\V}[2]{\mathbb{V}_{\mathbb{#1}#2}}

\DeclareMathOperator{\ball}{\mathfrak{B}}

\DeclareMathOperator{\ord}{Ord}

\DeclareMathOperator{\pow}{pow}

\numberwithin{equation}{section}

\begin{document}
	\begin{center}
		
		{\large {{\bf
					Class of extensions of real field and their topological properties.
		}}}
	\end{center}
	
	\begin{center}
		Alexandrov E. V.\\
		evgenyaexandrov66@gmail.com\\
	\end{center}
	
	{\centerline {\bf Abstract}}
	
	Proper classes of extensions of real field was defined and topological properties of these extensions were studied. These extensions can be connected, in this case such set is not closed under binary operations (addition and multiplication), and not connected, in this case this extension is linearly ordered field. In the future these constructions can be applied to building measure that <<feels>> set of zero Lebesgue measure.\\
	\hspace{1cm}
	
	\section{Introduction}
	
	In this paper we will construct extensions of sets $\mathbb{N}$, $\mathbb{Z}$, $\mathbb{Q}$ and $\mathbb{R}$. For each of these four sets a proper class of extensions, indexed by ordinals, will be constructed. Extensions of $\mathbb{N}$ and $\mathbb{Z}$ will contain infinitely large numbers, while extensions of $\mathbb{Q}$ and $\mathbb{R}$ will contain infinitely large numbers and infinitely small numbers.
	
	There are two types of extensions of $\mathbb{R}$:
	\begin{enumerate}
		\item Topologically connected extensions $\V{\overline{R}}{\alpha}$ ($\alpha$ -- ordinal index), that will not be closed under addition and hence will not be fields.
		
		\item Topologically not connected extensions $\V{R}{\alpha}$ ($\alpha$ -- ordinal index), that will be fields made from connected extensions by removing some points.
	\end{enumerate}
	
	These extensions were built using Von Neumann--Bernays--Gödel set of axioms. Nonstandard models were not used here.
	
	Also, we can denote $\V{R}{\infty}=\bigcup\limits_{\alpha\text{ -- ordinal}}\V{R}{\alpha}$, here $\V{R}{\alpha}$ -- extension of real set that is field. Author suggests, that $\V{R}{\infty}$ is isomorphic to class of surreal numbers, but this claim was not proven.
	
	\section{Some notation}
	In this document we will assume that $0\in\mathbb{N}$.

Infinite cardinals will be designated as $\aleph_{\alpha}$. Class of all ordinal numbers will be designated as $\ord$, $\omega_{\alpha}$ -- first ordinal number with cardinality $\aleph_{\alpha}$.


\begin{equation}
	\ord(\alpha)=\{\beta\in\ord|\beta<\alpha\}
\end{equation}
\begin{equation}
	\W{F}{\alpha}^+=\{x\in\W{F}{\alpha}|x>0\}
\end{equation}
similarly
\begin{equation}
	\V{F}{\alpha}^+=\{x\in\V{F}{\alpha}|x>0\}
\end{equation}
also
\begin{equation}
	\W{F}{\alpha}^{\geqslant 0}=\{x\in\W{F}{\alpha}|x\geqslant 0\}
\end{equation}
\begin{equation}
	\V{F}{\alpha}^{\geqslant 0}=\{x\in\V{F}{\alpha}|x\geqslant 0\}
\end{equation}

Let $f$ be some function, then $\ran(f)$ is range of $f$ and $\dom(f)$ is domain of $f$.

If $A$ and $B$ are sets and there is binary operation $*:A\times B\to C$ then
\begin{equation}
	(A*B)=\{x\in C|\exists y\in A\;\;\exists z\in B:x=y*z\}
\end{equation}
E.g. let $A,B\subset\mathbb{R}$ then
\begin{equation}
	(A+B)=\{x\in\mathbb{R}|\exists y\in A\;\;\exists z\in B:x=y+z\}
\end{equation}

If on set $S$ there is defined binary operation $*:S\times S\to S$, $A\subseteq S$ and $x\in S$, then
\begin{equation}
	\begin{cases}
		B=A*x\Leftrightarrow B=\{y\in S|\exists z\in A:y=z*x\}\\
		C=x*A\Leftrightarrow C=\{y\in S|\exists z\in A:y=x*z\}\\
	\end{cases}
\end{equation}
E.g. let $A\subset\mathbb{R}$ and $x\in\mathbb{R}$ then
\begin{equation}
	A+x=\{y\in\mathbb{R}|\exists z\in A:y=z+x\}
\end{equation}

If $f$ is some function and $A\subseteq\dom(f)$, then $f(A)$ is the following set:
\begin{equation}
	f(A)=\{x\in\ran(f)|\exists y\in A:f(y)=x\}
\end{equation}


	
	\section{$\W{N}{\alpha}$ and $\V{N}{\alpha}$}\label{NaturalNumbersSection}
	Everywhere in his document will be used axiomatics of Von Neumann–Bernays–Gödel.

Firstly we will consider class of ordinal numbers $\ord$ and binary operations on it.

We will define <<standard>> operations on $\ord$ following~\cite{GaisiTakeutiSetTheory}. Unlike other literature (including~\cite{GaisiTakeutiSetTheory}) we will denote <<standard>> addition, multiplication and exponentiation as $\oplus$, $\odot$ and $\pow$. Also here function that returns next ordinal we will denote as $S(\alpha)$.

\begin{definition}[Limit ordinal]
	Limit ordinal is ordinal number that has not got previous ordinal.
\end{definition}

\begin{definition}[<<Standard>> addition on $\ord$]
	We will define <<standard>> addition with transfinite recursion. Let $\alpha\in\ord$, then
	\begin{enumerate}
		\item $\alpha\oplus 0=\alpha$;
		\item $\alpha\oplus S(\beta)=S(\alpha\oplus\beta)$;
		\item if $\beta$ is limit ordinal, then $\alpha\oplus\beta=\sup\limits_{\gamma<\beta}(\alpha\oplus\gamma)$.
	\end{enumerate}
\end{definition}

\begin{definition}[<<Standard>> multiplication on $\ord$]
	We will define <<standard>> multiplication with transfinite recursion. Let $\alpha\in\ord$, then
	\begin{enumerate}
		\item $\alpha\odot 0=0$;
		\item $\alpha\odot S(\beta)=\alpha\odot\beta\oplus\alpha$;
		\item if $\beta$ is limit ordinal, then $\alpha\odot\beta=\sup\limits_{\gamma<\beta}(\alpha\odot\gamma)$.
	\end{enumerate}
\end{definition}

\begin{definition}[<<Standard>> exponentiation on $\ord$]\label{StandExpWNDef}
	We will define <<standard>> exponentiation with transfinite recursion. Let $\alpha\in\ord$, then
	\begin{enumerate}
		\item $\pow(\alpha,0)=1$;
		\item $\pow(\alpha,S(\beta))=\pow(\alpha,\beta)\odot\alpha$;
		\item if $\beta$ is limit ordinal, then $\pow(\alpha,\beta)=\sup\limits_{\gamma<\beta}(\pow(\alpha,\gamma))$
	\end{enumerate}
\end{definition}

According \cite{Kurat} and \cite{GaisiTakeutiSetTheory} these operations have following properties:
\begin{enumerate}
	\item addition is not commutative;
	\item addition is associative;
	\item multiplication is not commutative;
	\item multiplication is associative;
	\item for addition and multiplication there is left distributivity:
	\begin{equation}
		\alpha\odot(\beta\oplus\gamma)=\alpha\odot\beta\oplus\alpha\odot\gamma
	\end{equation}
\end{enumerate}

According to~\cite{Kurat} and~\cite{SuperJac} every ordinal number can be represented in unique way in Cantor normal form. Let $\alpha\in\ord$, then exists finite strictly decreasing sequence of ordinal numbers $\eta_k$ and sequence of natural numbers $a_k$ (both have same length $N\in\mathbb{N}$) such that

\begin{equation}
	\alpha=\pow(\omega_0,\eta_0)\odot a_0\oplus\ldots\oplus\pow(\omega_{N-1},\eta_{N-1})\odot a_{N-1}
\end{equation}

Following~\cite{Kurat} and~\cite{SuperJac} we can define natural addition and multiplication:
\begin{definition}[Natural addition]
	Let $\alpha,\beta\in\ord$ and in Cantor normal form (without loss of generality we can assume that lengths of these forms are equal and $\eta_k$ is the same for both numbers):
	\begin{equation}
		\begin{cases}
			\alpha=\bigoplus\limits_{k=0}^{N-1}\pow(\omega_0,\eta_k)\odot a_k\\
			\beta=\bigoplus\limits_{k=0}^{N-1}\pow(\omega_0,\eta_k)\odot b_k
		\end{cases}
	\end{equation}
	Then
	\begin{equation}
		\alpha+\beta=\bigoplus\limits_{k=0}^{N-1}\pow(\omega_0,\eta_k)\odot(a_k+b_k)
	\end{equation}
\end{definition}

According to~\cite{Kurat} and~\cite{SuperJac} natural addition is commutative and associative. Then we can notice that Cantor normal form can be rewritten as
\begin{equation}
	\alpha=\sum\limits_{k=0}^{N-1}\pow(\omega_0,\eta_k)\odot a_k
\end{equation}
Here $\sum$ is natural addition. In this form we can permute terms in arbitrary way.

\begin{definition}[Natural multiplication]
	Let $\alpha,\beta\in\ord$ and in Cantor normal form:
	\begin{equation}
		\begin{cases}
			\alpha=\sum\limits_{k=0}^{N-1}\pow(\omega_0,\eta_k)\odot a_k\\
			\beta=\sum\limits_{k=0}^{M-1}\pow(\omega_0,\gamma_k)\odot b_k
		\end{cases}
	\end{equation}
	Then
	\begin{equation}
		\alpha\cdot\beta=\sum\limits_{k=0}^{N-1}\sum\limits_{p=0}^{M-1}\pow(\omega_0,\eta_k+\gamma_p)\odot(a_k\cdot b_k)
	\end{equation}
\end{definition}

According to~\cite{Kurat} and~\cite{SuperJac} natural multiplication is commutative and associative. Also natural addition and multiplication are distributive.

Because <<standard>> multiplication of arbitrary ordinal to natural number equals to natural multiplication to same natural number, then we can rewrite Cantor normal form in the following manner:
\begin{equation}
	\alpha=\sum\limits_{k=0}^{N-1}\pow(\omega_0,\eta_k)\cdot a_k
\end{equation}

In~\cite{SuperJac} it is proven that natural exponentiation on ordinal numbers does not exist.

\begin{definition}[$\W{N}{\alpha}$]
	\begin{equation}
		\W{N}{\alpha}=\ord(\omega_{\alpha})
	\end{equation}
	addition and multiplication on $\W{N}{\alpha}$ defined as natural addition and multiplication.
\end{definition}

\begin{definition}[$\V{N}{\alpha}$]
	\begin{equation}
		\V{N}{\alpha}=\ord(\pow(\omega_0,\pow(\omega_0,\alpha)))
	\end{equation}
	addition and multiplication on $\V{N}{\alpha}$ defined as natural addition and multiplication.
\end{definition}

\begin{definition}[Principal numbers of addition]
	$\alpha\in\W{N}{\infty}$ is principal number of addition if
	\begin{equation}
		\begin{cases}
			\beta<\alpha\\
			\gamma<\alpha\\
		\end{cases}
		\Rightarrow
		\beta+\gamma<\alpha
	\end{equation}
\end{definition}

\begin{definition}[Principal number of multiplication]
	$\alpha\in\W{N}{\infty}$ is principal number of multiplication if
	\begin{equation}
		\begin{cases}
			\beta<\alpha\\
			\gamma<\alpha\\
		\end{cases}
		\Rightarrow
		\beta\cdot\gamma<\alpha
	\end{equation}
\end{definition}

\begin{theorem}[About principal numbers of addition]\label{PrincNumAdditionTheor}
	Let $\alpha\in\W{N}{\infty}$, then
	\begin{equation}
		\exists\xi\in\W{N}{\infty}:\alpha=\pow(\omega_0,\xi)
	\end{equation}
	if and only if
	\begin{equation}
		\forall\gamma<\alpha,\forall\beta<\alpha\;\;\gamma+\beta<\alpha
	\end{equation}
\end{theorem}

\textbf{Proof}\\

Let $\alpha=\pow(\omega_0,\delta)$ where $\delta\in\W{N}{\infty}$. Because $\beta<\alpha$ and $\gamma<\alpha$, it is obvious that $\beta$ and $\gamma$ have power of $\omega_0$ in Cantor normal form less then $\delta$, i.e.
\begin{equation}
	\begin{cases}
		\beta=\sum\limits_{k=0}^{N-1}\pow(\omega_0,\eta_k)\cdot b_k<\alpha\\
		\gamma=\sum\limits_{k=0}^{N-1}\pow(\omega_0,\eta_k)\cdot c_k<\alpha\\
	\end{cases}
	\Rightarrow
	\max\limits_{k<N}(\eta_k)<\delta
\end{equation}
Here some of $c_k$ or $b_k$ can be zero. But if we sum $\beta$ and $\gamma$ there will not be terms with $\pow(\omega_0,\xi)$ where $\xi\geqslant\delta$ in Cantor normal form of sum. Hence if $\alpha=\pow(\omega_0,\delta)$ then $\forall\beta<\alpha\;\;\forall\gamma<\alpha\;\;\beta+\gamma<\alpha$.

Now let $\alpha\in\W{N}{\infty}$ such that
\begin{equation}
	\begin{cases}
		\beta<\alpha\\
		\gamma<\alpha\\
	\end{cases}
	\Rightarrow
	\beta+\gamma<\alpha
\end{equation}
Let us assume that $\alpha$ has in Cantor normal form number of terms greater then one:
\begin{equation}
	\alpha=\sum\limits_{k=0}^{N-1}\pow(\omega_0,\eta_k)\cdot a_k
\end{equation}
Let $\eta_0=\max\limits_{k<N}(\eta_k)$, then consider
\begin{equation}
	\begin{cases}
		\beta=\pow(\omega_0,\eta_0)\cdot a_0\\
		\gamma=\beta=\pow(\omega_0,\eta_0)\cdot a_0\\
	\end{cases}
\end{equation}
In this case
\begin{equation}
	\begin{cases}
		\beta<\alpha\\
		\gamma<\alpha\\
		\beta+\gamma=\pow(\omega_0,\eta_0)\cdot 2\cdot a_0>\alpha\\
	\end{cases}
\end{equation}
Hence $N=1$, i.e.
\begin{equation}
	\alpha=\pow(\omega_0,\eta_0)\cdot a_0
\end{equation}
Let us now prove that $a_0=1$. We will use proof by contradiction, let $a_0>1$, then we can choose $\beta$ and $\gamma$ in the following manner
\begin{equation}
	\begin{cases}
		\beta=\pow(\omega_0,\eta_0)\cdot(a_0-1)\\
		\gamma=\pow(\omega_0,\eta_0)\cdot(a_0-1)+\pow(\omega_0,\eta_1)\\
	\end{cases}
\end{equation}
Here we mean that $\eta_1<\eta_0$. $a_0-1\in\mathbb{N}$ because $a_0>1$. Thus
\begin{equation}
	\beta+\gamma=\pow(\omega_0,\eta_0)\cdot(2a_0-2)+\pow(\omega_0,\eta_1)
\end{equation}
Because $a_0>1$, then $2a_0-2\geqslant a_0$, hence $\pow(\omega_0,\eta_0)(2a_0-2)\geqslant\alpha$. Because $\pow(\omega_0,\eta_1)>0$, then
\begin{equation*}
	\pow(\omega_0,\eta_0)(2a_0-2)+\pow(\omega_0,\eta_1)>\alpha
	\Rightarrow
\end{equation*}
\begin{equation*}
	\Rightarrow
	\beta+\gamma>\alpha
\end{equation*}
Hence we have contradiction and hence $a_0=1$, i.e. $\alpha=\pow(\omega_0,\eta_0)$.

$\square$

\begin{theorem}[About principal numbers of multiplication]\label{PrincNumMultWNTheor}
	Let $\alpha\in\W{N}{\infty}$, then
	\begin{equation}
		\exists\xi\in\W{N}{\infty}:\alpha=\pow(\omega_0,\pow(\omega_0,\xi))
	\end{equation}
	if and only if
	\begin{equation}
		\forall\gamma<\alpha,\forall\beta<\alpha\;\;\beta\cdot\gamma<\alpha
	\end{equation}
\end{theorem}

\textbf{Proof}\\

Let $\alpha=\pow(\omega_0,\pow(\omega_0,\delta))$ where $\delta\in\W{N}{\infty}$. Because $\beta<\alpha$ and $\gamma<\alpha$, it is obvious that $\beta$ and $\gamma$ have power of $\omega_0$ in Cantor normal form less then $\pow(\omega_0,\delta)$, i.e.
\begin{equation}
	\begin{cases}
		\beta=\sum\limits_{k=0}^{N-1}\pow(\omega_0,\eta_k)\cdot b_k<\alpha\\
		\gamma=\sum\limits_{k=0}^{M-1}\pow(\omega_0,\psi_k)\cdot c_k<\alpha\\
	\end{cases}
	\Rightarrow
	\begin{cases}
		\max\limits_{k<N}(\eta_k)=\eta_0<\pow(\omega_0,\delta)\\
		\max\limits_{k<M}(\psi_k)=\psi_0<\pow(\omega_0,\delta)
	\end{cases}
\end{equation}
Here every $b_k$ and $c_k$ does not equal to zero. Hence
\begin{equation}
	\gamma\cdot\beta=\sum\limits_{k=0}^{N-1}\sum\limits_{p=0}^{M-1}\pow(\omega_0,\eta_k+\psi_p)\cdot(b_k\cdot c_p)
\end{equation}
Thus
\begin{equation}
	\max\limits_{k<N\wedge p<M}(\eta_k+\psi_p)=\eta_0+\psi_0
\end{equation}
Because $\pow(\omega_0,\delta)$ -- principal number of addition, then $\eta_0+\psi_0<\pow(\omega_0,\delta)$, hence
\begin{equation}
	\beta\cdot\gamma<\pow(\omega_0,\pow(\omega_0,\delta))=\alpha
\end{equation}

Now let $\alpha\in\W{N}{\infty}$ such that
\begin{equation}
	\begin{cases}
		\beta<\alpha\\
		\gamma<\alpha\\
	\end{cases}
	\Rightarrow
	\beta\cdot\gamma<\alpha
\end{equation}
Consider Cantor normal form of $\alpha$:
\begin{equation}
	\alpha=\sum\limits_{k=0}^{N-1}\pow(\omega_0,\eta_k)\cdot a_k
\end{equation}
Let us prove that $N=1$ using proof by contradiction. Let us assume that $N>1$. Let $\eta_0=\max\limits_{k<N}(\eta_k)$, then consider
\begin{equation}
	\begin{cases}
		\beta=\pow(\omega_0,\eta_0)\cdot a_0\\
		\gamma=\beta=\pow(\omega_0,\eta_0)\cdot a_0\\
	\end{cases}
\end{equation}
In this case
\begin{equation}
	\begin{cases}
		\beta<\alpha\\
		\gamma<\alpha\\
		\beta\cdot\gamma=\pow(\omega_0,\eta_0+\eta_0)\cdot a_0^2>\alpha\\
	\end{cases}
\end{equation}
$\pow(\omega_0,\eta_0+\eta_0)\cdot a_0^2>\alpha$ because $\eta_0+\eta_0>\eta_0$. Thus we have contradiction and hence $N=1$. I.e. $\alpha=\pow(\omega_0,\eta_0)\cdot a_0$.

Let us prove that $\eta_0=\pow(\omega_0,\delta)$ for some $\delta$ using proof by contradiction. Let us assume that $\eta_0$ is not principal number of addition, then consider
\begin{equation}
	\begin{cases}
		\beta=\pow(\omega_0,\varphi_0)\\
		\gamma=\pow(\omega_0,\psi_0)\\
		\varphi_0<\eta_0\\
		\psi_0<\eta_0\\
		\varphi_0+\psi_0>\eta_0\\
	\end{cases}
\end{equation}
Obviously
\begin{equation}
	\gamma\cdot\beta=\pow(\omega_0,\varphi_0+\psi_0)
\end{equation}
Because $\varphi_0+\psi_0>\eta_0$, then $\gamma\cdot\beta>\alpha$, hence we have contradiction and hence $\eta_0$ is principal number of addition, i.e. $\eta_0=\pow(\omega_0,\delta)$. Thus we have proven that $\alpha=\pow(\omega_0,\pow(\omega_0,\delta))\cdot a_0$.

Let us prove that $a_0=1$ using proof by contradiction, let us assume that $a_0>1$, then consider
\begin{equation}
	\begin{cases}
		\beta=\pow(\omega_0,\pow(\omega_0,\delta))\cdot(a_0-1)\\
		\gamma=\beta=\pow(\omega_0,\pow(\omega_0,\delta))\cdot(a_0-1)\\
	\end{cases}
\end{equation}
Then
\begin{equation}
	\beta\cdot\gamma=\pow(\omega_0,2\cdot\pow(\omega_0,\delta))\cdot(a_0-1)^2
\end{equation}
Because $2\cdot\pow(\omega_0,\delta)>\pow(\omega_0,\delta)$, then $\beta\cdot\gamma>\alpha$, hence we have contradiction and hence $a_0=1$, i.e. $\alpha=\pow(\omega_0,\pow(\omega_0,\delta))$.

$\square$

Let us notice that any principal number of multiplication is principal number of addition, reverse is not true.

Here we can see, because $\omega_{\alpha}$ and $\pow(\omega_0,\pow(\omega_0,\alpha))$ are both principal numbers of multiplication, than for all $\alpha\in\ord$ both $\W{N}{\alpha}$ and $\V{N}{\alpha}$ are closed under addition and multiplication.

Also let us notice that $\W{N}{\alpha}=\V{N}{\omega_{\alpha}}$.
	
	\section{$\W{Z}{\alpha}$ and $\V{Z}{\alpha}$}
	In this section will be defined two classes of extensions of integers: $\W{Z}{\alpha}$ and $\V{Z}{\alpha}$, where $\alpha\in\ord$. Every extension $\W{Z}{\alpha}$ will be defined on a basis of $\W{N}{\alpha}$ ($\V{Z}{\alpha}$ -- similarly). When moving from $\W{N}{\alpha}$ to $\W{Z}{\alpha}$ not only negative integers will appear (e.g. $\omega_0-1>0$ but $(\omega_0-1)\notin\W{N}{\alpha}$).

\subsection{Definition}

\begin{definition}[$Z_{\alpha}$]
	$Z_{\alpha}$ -- set of ordered pairs of elements of $\V{N}{\alpha}$, that is
	\begin{equation}
		a,b\in\V{N}{\alpha}\Rightarrow(a,b)\in Z_{\alpha}
	\end{equation}
\end{definition}

Let us define an equivalence relation on $Z_{\alpha}$, let
\begin{equation}
	(a,b)\in Z_{\alpha}\wedge(c,d)\in Z_{\alpha}
\end{equation}
then
\begin{equation}
	(a,b)\sim(c,d)\Leftrightarrow a+d=b+c
\end{equation}

\begin{theorem}[About equivalence on $Z_{\alpha}$]
	Relation just defined is equivalence relation.
\end{theorem}

\textbf{Proof}\\

\begin{enumerate}
	\item $(a,b)\sim(a,b)$ follows from $a+b=b+a$;
	\item $(a,b)\sim(c,d)\Rightarrow(c,d)\sim(a,b)$ -- obvious;
	\item $(a,b)\sim(c,d)\wedge(c,d)\sim(e,g)\Rightarrow(a,b)\sim(e,g)$:\\
	According definition
	\begin{equation}
		\begin{cases}
			a+d=b+c\\
			c+g=d+e\\
		\end{cases}
	\end{equation}
	let us sum these equations:
	\begin{equation}
		a+g+d+c=b+e+d+c\Rightarrow a+g=b+e\Rightarrow(a,b)\sim(e,g)
	\end{equation}
\end{enumerate}

$\square$

Let us define addition and multiplication on $Z_{\alpha}$:
\begin{equation}
	(a,b),(c,d)\in Z_{\alpha}\Rightarrow
	\begin{cases}
		(a,b)+(c,d)=(a+c,b+d)\\
		(a,b)\cdot(c,d)=(a\cdot c+b\cdot d,a\cdot d+b\cdot c)\\
	\end{cases}
\end{equation}
Commutativity and associativity are obvious.

\begin{theorem}[About consistency of operations on $Z_{\alpha}$ with equivalence]
	\begin{equation}
		\begin{cases}
			(a_1,b_1)\sim(a_2,b_2)\\
			(c_1,d_1)\sim(c_2,d_2)\\
		\end{cases}
		\Rightarrow
		\begin{cases}
			(a_1,b_1)+(c_1,d_1)\sim(a_2,b_2)+(c_2,d_2)\\
			(a_1,b_1)\cdot(c_1,d_1)\sim(a_2,b_2)\cdot(c_2,d_2)\\
		\end{cases}
	\end{equation}
\end{theorem}

\textbf{Proof}\\

Let us proove consistency of addition. According to definition
\begin{equation}
	\begin{cases}
		a_1+b_2=a_2+b_1\\
		c_1+d_2=c_2+d_1\\
	\end{cases}
	\Rightarrow
\end{equation}
\begin{equation}
	\Rightarrow
	(a_1+c_1)+(b_2+d_2)=(a_2+c_2)+(b_1+d_1)
	\Rightarrow
\end{equation}
\begin{equation}
	\Rightarrow
	(a_1,b_1)+(c_1,d_1)\sim(a_2,b_2)+(c_2,d_2)
\end{equation}

Proof of consistency of multiplication is similar.

$\square$

\begin{definition}[$\V{Z}{\alpha}$]
	$\V{Z}{\alpha}=Z_{\alpha}/_{\sim}$, factorization is based on just defined equivalence.
\end{definition}

Because addition and multiplication are consistent with equivalence, then they are correctly defined on $\V{Z}{\alpha}$.

\begin{definition}[$\W{Z}{\alpha}$]
	\begin{equation*}
		\W{Z}{\alpha}=\V{Z}{\omega_{\alpha}}
	\end{equation*}
\end{definition}

\subsection{Order definition}

Let us define order relation on $\V{Z}{\alpha}$.

\begin{definition}[Order relation on $\V{Z}{\alpha}$]
	\begin{equation}
		A,B\in\V{Z}{\alpha}
	\end{equation}
	\begin{equation}
		A\geqslant B\Leftrightarrow
		\forall(a_1,a_2)\in A\;\;\forall(b_1,b_2)\in B\;\;a_1+b_2\geqslant a_2+b_1
	\end{equation}
\end{definition}

\begin{theorem}[About linear order on $\V{Z}{\alpha}$]
	Just defined order on $\V{Z}{\alpha}$ is linear order.
\end{theorem}

\textbf{Proof}\\

\begin{enumerate}
	\item reflexivity is obvious;
	\item $A\geqslant B\wedge B\geqslant C\Rightarrow$
	\begin{equation}
		\forall(a_1,a_2)\in A\;\;\forall(b_1,b_2)\in B\;\;\forall(c_1,c_2)\in C
		\begin{cases}
			a_1+b_2\geqslant a_2+b_1\\
			b_1+c_2\geqslant b_2+c_1\\
		\end{cases}
		\Rightarrow
	\end{equation}
	\begin{equation}
		\Rightarrow
		a_1+c_2+b_1+b_2\geqslant a_2+c_1+b_1+b_2
		\Rightarrow
	\end{equation}
	\begin{equation}
		\Rightarrow
		a_1+c_2\geqslant a_2+c_1\Rightarrow A\geqslant C
	\end{equation}
	\item $A\geqslant B\wedge B\geqslant A\Rightarrow$
	\begin{equation}
		\Rightarrow
		\forall(a_1,a_2)\in A\;\;\forall(b_1,b_2)\in B
		\begin{cases}
			a_1+b_2\geqslant a_2+b_1\\
			a_2+b_1\geqslant a_1+b_2\\
		\end{cases}
		\Rightarrow
	\end{equation}
	\begin{equation}
		\Rightarrow
		a_1+b_2=a_2+b_1
		\Rightarrow
		(a_1,a_2)\sim(b_1,b_2)\Rightarrow A=B
	\end{equation}
	\item Let $A,B\in\W{Z}{\alpha}$, $(a_1,a_2)\in A$, $(b_1,b_2)\in B$ and
	\begin{equation}
		a_1+b_2\geqslant a_2+b_1
	\end{equation}
	Consider
	\begin{equation}
		(a'_1,a'_2)\in A\Leftrightarrow
		\exists f\in\W{N}{\alpha}:
		\begin{cases}
			a'_1=a_1+f\\
			a'_2=a_2+f\\
		\end{cases}
	\end{equation}
	thus, comparability of any elements of $\V{Z}{\alpha}$ follows from fact that
	\begin{equation}
		a_1+b_2\geqslant a_2+b_1\Leftrightarrow
		(a_1+f)+b_2\geqslant (a_2+f)+b_1
	\end{equation}
\end{enumerate}

$\square$

Let us note, that after transition from $\V{N}{\alpha}$ to $\V{Z}{\alpha}$ not only negative elements appeared. For example:
\begin{equation}
	\omega_0-1>0
\end{equation}
but
\begin{equation}
	(\omega_0-1)\notin\V{N}{\alpha}
\end{equation}

Because $\W{Z}{\alpha}=\V{Z}{\omega_{\alpha}}$, then order is defined on $\W{Z}{\alpha}$ and has same properties.

\subsection{Embedding}

\begin{definition}[Embedding $\V{N}{\alpha}\subset\V{Z}{\alpha}$]
	Let $A\in \V{N}{\alpha}$ and $B\in\V{Z}{\alpha}$, then
	\begin{equation}
		A=B\Leftrightarrow(A,0)\in B
	\end{equation}
\end{definition}

\begin{definition}[Embedding $\V{Z}{\alpha}\subset\V{Z}{\beta}$]
	Let $\alpha<\beta$, $A\in\V{Z}{\alpha}$ and $B\in\V{Z}{\beta}$, then
	\begin{equation}
		A=B\Leftrightarrow
		\exists(p,q)\in A\wedge\exists(g,f)\in B:(p,q)=(g,f)
	\end{equation}
\end{definition}

	\section{$\W{Q}{\alpha}$ and $\V{Q}{\alpha}$}
	In this section will be defined two classes of extensions of rational numbers $\W{Q}{\alpha}$ and $\V{Q}{\alpha}$ where $\alpha\in\ord$. Every extension $\V{Q}{\alpha}$ will be defined on a basis of $\V{Z}{\alpha}$.

\subsection{Definition}

\begin{definition}[$Q_{\alpha}$]
	$Q_{\alpha}$ -- set of ordered pairs $(a,b)$ where $a\in\V{Z}{\alpha}$ and $b\in\V{Z}{\alpha}\setminus\{0\}$.
\end{definition}

Let us define binary operations on $Q_{\alpha}$:
\begin{equation}
	(a,b)\in Q_{\alpha}\wedge(c,d)\in Q_{\alpha}
\end{equation}
\begin{equation}
	(a,b)\cdot(c,d)=(a\cdot c,b\cdot d)
\end{equation}
\begin{equation}
	(a,b)+(c,d)=(a\cdot d+b\cdot c,b\cdot d)
\end{equation}
\begin{equation}
	(a,b)-(c,d)=(a,b)+(-c,d)
\end{equation}

Commutativity, associativity and distributivity of addition and multiplication are obvious.

Let us define equivalence relation on $Q_{\alpha}$:
\begin{definition}[Equivalence on $Q_{\alpha}$]
	\begin{equation}
		(a,b),(c,d)\in Q_{\alpha}
	\end{equation}
	\begin{equation}
		(a,b)\sim(c,d)\Leftrightarrow
		a\cdot d=b\cdot c
	\end{equation}
\end{definition}

\begin{theorem}[About equivalence on $Q_{\alpha}$]
	Just defined relation is equivalence relation.
\end{theorem}

\textbf{Proof}\\

Reflexivity and symmetry are obvious. Let us check transitivity.
\begin{equation}
	\begin{cases}
		(a,b)\sim(c,d)\\
		(c,d)\sim(e,f)\\
	\end{cases}
	\Rightarrow
	\begin{cases}
		a\cdot d=b\cdot c\\
		c\cdot f=d\cdot e\\
	\end{cases}
\end{equation}
Let us multiply these equations:
\begin{equation}
	a\cdot d\cdot c\cdot f=b\cdot c\cdot d\cdot e
	\Rightarrow
\end{equation}
\begin{equation}
	a\cdot f=b\cdot e\Rightarrow
	(a,b)\sim(e,f)
\end{equation}

$\square$

Consistency of operations with equivalence is obvious.

\begin{definition}[$\V{Q}{\alpha}$]
	\begin{equation}
		\V{Q}{\alpha}=Q_{\alpha}/_{\sim}
	\end{equation}
\end{definition}

\begin{definition}[$\W{Q}{\alpha}$]
	\begin{equation}
		\W{Q}{\alpha}=\V{Q}{\omega_{\alpha}}
	\end{equation}
\end{definition}

Because addition and multiplication are consistent with equivalence, then they are correctly defined on $\V{Q}{\alpha}$.

\subsection{Order definition}

\begin{definition}[Order on $\V{Q}{\alpha}$]
	\begin{equation}
		A,B\in\V{Q}{\alpha}
	\end{equation}
	\begin{enumerate}
		\item
		\begin{equation}
			A\geqslant 0\Leftrightarrow
			\forall(a,b)\in A\;\;a\cdot b\geqslant 0
		\end{equation}
		\item
		\begin{equation}
			A\geqslant B\Leftrightarrow
			A-B\geqslant 0
		\end{equation}
	\end{enumerate}
\end{definition}

\begin{theorem}[About linear order on $\V{Q}{\alpha}$]
	Just define order is linear order on $\V{Q}{\alpha}$.
\end{theorem}

\textbf{Proof}

\begin{enumerate}
	\item reflexivity is obvious;
	\item $A\geqslant B\wedge B\geqslant C\Rightarrow$
	\begin{equation}
		\Rightarrow
		\begin{cases}
			A-B\geqslant 0\\
			B-C\geqslant 0\\
		\end{cases}
		\Rightarrow
		A-B+B-C\geqslant 0\Rightarrow
	\end{equation}
	\begin{equation}
		\Rightarrow
		A-C\geqslant 0\Rightarrow A\geqslant C
	\end{equation}
	\item comparability of each elements follows from comparability of each element with zero;
	\item antisymmetry follows from comparability of each elements.
\end{enumerate}

$\square$

\subsection{$\W{Q}{\alpha}$ and $\V{Q}{\alpha}$ properties}

\begin{theorem}[About field of rational numbers]
	$\V{Q}{\alpha}$ is lineary ordered field.
\end{theorem}

\textbf{Proof}

Commutativity associativity and distributivity of addition and multiplication have already mentioned.

\begin{enumerate}
	\item Neutral element of addition is $\frac{0}{1}$:
	\begin{equation}
		(a,b)+(0,1)=(a\cdot 1+b\cdot 0,b\cdot 1)=(a,b)
	\end{equation}
	
	\item Inverse element of addition of $\frac{a}{b}$ is $\frac{-a}{b}$:
	\begin{equation}
		(a,b)+(-a,b)=(a\cdot b-a\cdot b,b^2)=(0,b^2)\sim(0,1)
	\end{equation}
	
	\item Neutral element of multiplication is $\frac{1}{1}$:
	\begin{equation}
		(a,b)\cdot(1,1)=(a\cdot 1,b\cdot 1)=(a,b)
	\end{equation}
	
	\item Inverse element of multiplication of $\frac{a}{b}$ ($a,b\ne 0$) is $\frac{b}{a}$:
	\begin{equation}
		(a,b)\cdot(b,a)=(a\cdot b,b\cdot a)\sim(1,1)
	\end{equation}
	
	\item Let $a,b,c\in\V{Q}{\alpha}$ and $a\geqslant b$. Consider
	\begin{equation}
		0\leqslant a-b=a+c-b+c\Rightarrow
	\end{equation}
	\begin{equation}
		\Rightarrow
		a+c\geqslant b+c
	\end{equation}
	
	\item claim that
	\begin{equation}
		\begin{cases}
			a\geqslant 0\\
			b\geqslant 0\\
		\end{cases}
		\Rightarrow
		a\cdot b\geqslant 0
	\end{equation}
	is obvious.
	
	\item Let $a,b,c\in\V{Q}{\alpha}^+$ and $a\geqslant b$. Consider
	\begin{equation}
		0\leqslant a-b=\frac{1}{c}\cdot c\cdot(a-b)
	\end{equation}
	because
	\begin{equation}
		\begin{cases}
			\frac{1}{c}>0\\
			c>0\\
			a-b\geqslant 0\\
		\end{cases}
	\end{equation}
	according to previous point we have
	\begin{equation}
		c\cdot(a-b)\geqslant 0
		\Rightarrow
	\end{equation}
	\begin{equation}
		\Rightarrow
		c\cdot a-c\cdot b\geqslant 0
		\Rightarrow
	\end{equation}
	\begin{equation}
		\Rightarrow
		a\cdot c\geqslant b\cdot c
	\end{equation}
\end{enumerate}

$\square$

\subsection{Embedding}

\begin{definition}[Embedding $\V{Z}{\alpha}\subset\V{Q}{\alpha}$]
	Let $A\in\V{Z}{\alpha}$ and $B\in\V{Q}{\alpha}$, then
	\begin{equation}
		A=B\Leftrightarrow B=(A,1)
	\end{equation}
\end{definition}

\begin{definition}[Embedding $\V{Q}{\alpha}\subset\V{Q}{\beta}$]
	Let $\alpha<\beta$, $A\in\V{Q}{\alpha}$ and $B\in\V{Q}{\beta}$ then
	\begin{equation}
		A=B\Leftrightarrow
		\exists(p,q)\in A\wedge\exists(g,f)\in B:(p,q)=(g,f)
	\end{equation}
\end{definition}

	\section{$\W{\overline{R}}{\alpha}$ and $\V{\overline{R}}{\alpha}$}\label{ExtRealDefSection}
	In this section will be defined two classes of extensions of set of real numbers -- $\W{\overline{R}}{\alpha}$ and $\V{\overline{R}}{\alpha}$ where $\alpha\in\ord$. About sets $\W{\overline{R}}{\alpha}$ and $\V{\overline{R}}{\alpha}$ you may think as about analogue of extended real line $\overline{\mathbb{R}}=[-\infty,+\infty]$, what is not field too.

\subsection{Definition, order}

\begin{definition}[Dedeking cut]
	Let $A,B\subset Q$ where $Q$ is linerly ordered set, then $(A,B)$ is Dedekind cut on $Q$ if next conditions are met:
	\begin{enumerate}
		\item
		\begin{equation}
			\forall a\in A\;\;\forall b\in B\;\;a<b
		\end{equation}
		\item set
		\begin{equation}
			Q\setminus\left(A\cup B\right)
		\end{equation}
		consists of one or zero elements
		\item
		\begin{equation}
			\begin{cases}
				\inf(B)\notin B\\
				\sup(A)\notin A\\
			\end{cases}
		\end{equation}
		\item
		\begin{equation}
			A\cap B=\varnothing
		\end{equation}
	\end{enumerate}
\end{definition}

\begin{definition}[Order on Dedekind cuts]
	Let $(A_1,B_1)$ and $(A_2,B_2)$ are Dedekind cuts on $Q$, then
	\begin{equation}
		(A_1,B_1)\geqslant(A_2,B_2)\Leftrightarrow A_2\subseteq A_1
	\end{equation}
\end{definition}

\begin{theorem}[About linear order on Dedeking cuts]
	Just defined order is linear order.
\end{theorem}

\textbf{Proof}\\

\begin{enumerate}
	\item reflexivity follows from this: $\forall A\;\;A\subseteq A$;
	\item transitivity follows from:
	\begin{equation}
		\begin{cases}
			A\subseteq B\\
			B\subseteq C\\
		\end{cases}
		\Rightarrow
		A\subseteq C
	\end{equation}
	\item antisymmetry follows from:
	\begin{equation}
		\begin{cases}
			A\subseteq B\\
			B\subseteq A\\
		\end{cases}
		\Rightarrow
		A=B
	\end{equation}
	\item comparability of each cuts $(A_1,B_1)$ and $(A_2,B_2)$ follows from fact that either $A_1\subseteq A_2$ or the $A_2\subseteq A_1$ or both is possible.
\end{enumerate}

$\square$

\begin{definition}[Improper extension of real field -- $\V{\overline{R}}{\alpha}$]
	If $D$ is set of Dedekind cuts on $\V{Q}{\alpha}$, then $\V{\overline{R}}{\alpha}$ is set of Dedekind cuts on $\V{Q}{\alpha}$.
	
	Particularly, if $A=\varnothing$, then $(A,B)=-\infty$, if $B=\varnothing$, then $(A,B)=+\infty$.
\end{definition}

\begin{definition}[Improper extension of real field -- $\W{\overline{R}}{\alpha}$]
	\begin{equation}
		\W{\overline{R}}{\alpha}=\V{\overline{R}}{\omega_{\alpha}}
	\end{equation}
\end{definition}

\subsection{Embedding}

\begin{definition}[Embedding $\V{Q}{\alpha}\subset\V{\overline{R}}{\alpha}$]\label{QREmbedDef}
	Let $x\in\V{Q}{\alpha}$ and $(A,B)\in\V{\overline{R}}{\alpha}$, then
	\begin{equation}
		x=(A,B)\Leftrightarrow x=\sup(A)
	\end{equation}
\end{definition}

\begin{definition}[Embedding $\V{\overline{R}}{\alpha}\subset\V{\overline{R}}{\beta}$]
	Let $\alpha<\beta$, let us define embedding function $f:\V{\overline{R}}{\alpha}\to\V{\overline{R}}{\beta}$ in the following way, let $x\in\V{\overline{R}}{\alpha}$:
	\begin{equation}
		f(x)=(C,D)\Leftrightarrow
		\begin{cases}
			\forall y\in C\;\;y<x\\
			\forall z\in D\;\;z>x\\
		\end{cases}
	\end{equation}
\end{definition}

\begin{theorem}[About embedding of $\V{\overline{R}}{\alpha}$]
	Just defined embedding function is single-valued function.
\end{theorem}

\textbf{Proof}\\

Let us use proof by contradiction. Let us assume that for some $(C,D)\in\V{\overline{R}}{\beta}$ exist $(A_1,B_1),(A_2,B_2)\in\V{\overline{R}}{\alpha}$ such that
\begin{enumerate}
	\item $(A_1,B_1)\ne(A_2,B_2)$;
	\item
	\begin{equation}
		\begin{cases}
			\forall y\in C\;\;y<(A_1,B_1)\wedge y<(A_2,B_2)\\
			\forall z\in D\;\;z>(A_1,B_1)\wedge z>(A_2,B_2)\\
		\end{cases}
	\end{equation}
\end{enumerate}

Without loss of generality we can assume that $(A_2,B_2)>(A_1,B_1)$.

Consider $S=A_2\cap B_1$, let $x\in S$. Because exists $\delta\in\V{Q}{\beta}^+$ such that
\begin{equation}
	\forall p\in\V{Q}{\alpha}\setminus\{0\}\;\;|\delta|<|p|
\end{equation}
then $y=x+\delta$ is located between $C$ and $D$. But in this case $(C,D)$ is not Dedekind cut, hence we have contradiction and hence embedding map is single-valued.

$\square$

\subsection{Topology}\label{SubsectionTopologyOfR}

\begin{definition}[Ball on $\V{\overline{R}}{\alpha}$]\label{BallOnExtRDef}
	$A\subset\V{\overline{R}}{\alpha}$ is ball on $\V{\overline{R}}{\alpha}$ if
	\begin{equation}
		\exists x,y\in\V{\overline{R}}{\alpha}:
		A=\{z\in\V{\overline{R}}{\alpha}|x<z<y\}
	\end{equation}
\end{definition}

\begin{theorem}[About base of topology on $\V{\overline{R}}{\alpha}$]
	Set of balls on $\V{\overline{R}}{\alpha}$ is topology base.
\end{theorem}

\textbf{Proof}\\

Let $B$ is set of balls on $\V{\overline{R}}{\alpha}$. It is anough to proove \cite{KelTop} that
\begin{equation}
	\forall U,V\in B\;\;\forall x\in U\cap V\;\;\exists W\in B:
	\begin{cases}
		x\in W\\
		W\subset U\cap V
	\end{cases}
\end{equation}
Indeed, consider $U,V\in B$ and $U\cap V\ne\varnothing$, by definition:
\begin{equation}
	\exists x_u,y_u,x_v,y_v\in\V{\overline{R}}{\alpha}:
	\begin{cases}
		U=\{z\in\V{\overline{R}}{\alpha}|x_u<z<y_u\}\\
		V=\{z\in\V{\overline{R}}{\alpha}|x_v<z<y_v\}\\
	\end{cases}
\end{equation}
then it is obvious that
\begin{equation}
	U\cap V=\{z\in\V{\overline{R}}{\alpha}|\max(x_u,x_v)<z<\min(y_u,y_v)\}\in B
\end{equation}

$\square$

\begin{theorem}[About infimum and supremum on $\V{\overline{R}}{\alpha}$]\label{AboutInfAndSupOnVRTheor}
	Each bounded set in $\V{\overline{R}}{\alpha}$ has supremum and infimum which belong $\V{\overline{R}}{\alpha}$.
\end{theorem}

\textbf{Proof}\\

Let $S=\subset\V{\overline{R}}{\alpha}$ be upper bounded set. Let us build following sets:
\begin{equation}
	\begin{cases}
		A=\{x\in\V{Q}{\alpha}|\exists p\in S:p>x\}\\
		B=\V{Q}{\alpha}\setminus A\\
	\end{cases}
\end{equation}
If exists $y=\min(B)$, then $\sup(S)=y$, otherwise $(A,B)$ is Dedekind cut and hence $(A,B)=\sup(S)$.

Proof for infimum is similar.

$\square$

\begin{theorem}[About connectedness of set of Dedeking cuts]
	$\V{\overline{R}}{\alpha}$ is connected set in order topology (just defined).
\end{theorem}

\textbf{Proof}\\

Exist two opened sets on $\V{\overline{R}}{\alpha}$ with empty boundary -- $\varnothing$ and $\V{\overline{R}}{\alpha}$, other sets have non-empty boundary, it follows from theorem \ref{AboutInfAndSupOnVRTheor}. Hence every opened set except $\varnothing$ and $\V{\overline{R}}{\alpha}$ do not coincides with its closure. Thus sets that are opened and closed simultaneously are only $\varnothing$ and $\V{\overline{R}}{\alpha}$, hence $\V{\overline{R}}{\alpha}$ is connected.

$\square$

\begin{theorem}[About Hausdorff topology on $\V{\overline{R}}{\alpha}$]
	Topology on $\V{\overline{R}}{\alpha}$ is Hausdorff topology.
\end{theorem}

\textbf{Proof}\\

Consider Dedekind cuts $x,y\in\V{\overline{R}}{\alpha}$: $x=(A,B)$, $y=(C,D)$, $x<y$. Because $x<y$ then $B\cap C\ne\varnothing$. Let us choose one point $z\in B\cap C$ and two points $z_0\in A$, $z_1\in D$. Consider balls
\begin{equation}
	\begin{cases}
		\ball_1=\{a\in\V{\overline{R}}{\alpha}|z_0<a<z\}\\
		\ball_2=\{b\in\V{\overline{R}}{\alpha}|z<b<z_1\}\\
	\end{cases}
\end{equation}
Obviously $\ball_1$ and $\ball_2$ are neighborhoods of $x$ and $y$ respectively. Also $\ball_1\cap\ball_2=\varnothing$. Thus every pair of different points have non-intersecting neighborhoods.

$\square$

\begin{definition}[Ball on linearly ordered field]
	Let $Q$ be some linearly ordered field, $x\in Q$ and $r\in Q^+$. then set
	\begin{equation}
		\ball(x,r)=\{y\in Q|r>|x-y|\}
	\end{equation}
	is a ball on $Q$.
\end{definition}

\begin{definition}[Improper elements]
	Let $(A,B)$ be Dedekind cut on $\V{Q}{\alpha}$ and $Q$ is set of balls on $\V{Q}{\alpha}$ such that
	\begin{equation}
		\ball\in Q\Leftrightarrow
		\begin{cases}
			A\cap\ball\ne\varnothing\\
			B\cap\ball\ne\varnothing\\
		\end{cases}
	\end{equation}
	Then $(A,B)$ is improper if
	\begin{equation}
		\exists d\in\V{Q}{\alpha}^+:\forall\ball\in Q\;\;\sup(\ball)-\inf(\ball)>d
	\end{equation}
	It means that infimum of ball diameters is greater then zero.
\end{definition}

\begin{theorem}[About existence of improper elements]
	There are improper elements in $\V{\overline{R}}{\alpha}$ if $\alpha>0$.
\end{theorem}

\textbf{Proof}\\

Let us give an example. Let $S$ be a set such that (here alpha can take arbitrary value greater then zero)
\begin{equation}
	S=\{x\in\V{Q}{\alpha}|\forall y\in\mathbb{Q}^+|x|<y\}
\end{equation}
Also let $(A,B)\in\V{\overline{R}}{\alpha}$ be Dedekind cut such that
\begin{equation}
	\begin{cases}
		A=\{x\in\V{Q}{\alpha}|x<\sup(S)\}\\
		B=\{x\in\V{Q}{\alpha}|x>\sup(S)\}\\
	\end{cases}
\end{equation}
Consider ball $\ball$ such that
\begin{equation}
	\begin{cases}
		A\cap\ball\ne\varnothing\\
		B\cap\ball\ne\varnothing\\
	\end{cases}
\end{equation}

Let us notice that
\begin{equation}
	\begin{cases}
		\inf(\ball)\in S\\
		\forall q\in S\;\;\sup(\ball)>q\\
	\end{cases}
\end{equation}
hence
\begin{equation}
	\forall q\in S\;\;\sup(\ball)-\inf(\ball)>q
\end{equation}
Claim follows from fact that $S$ contains elements greater then zero.

$\square$

\subsection{Binary operations}

\begin{definition}[Sum on Dedekind cuts]
	Let $(A,B)$ and $(C,D)$ be Dedekind cuts on $\V{Q}{\alpha}$, then pair $(E,F)$ is sum of $(A,B)$ and $(C,D)$ if
	\begin{equation}
		\begin{cases}
			E=A+C\\
			F=B+D\\
		\end{cases}
	\end{equation}
	here $A+C$ and $B+D$ are sums defined in section <<Some notation>>.
\end{definition}

\textbf{Properties of addition on Dedekind cuts}\\

\begin{enumerate}
	\item Sum of Dedekind cuts is not necessary a Dedekind cut.
	
	\item Sum on Dedekind cuts is commutative and associative (it follows from corresponding properties of $\V{Q}{\alpha}$).
	
	\item Exists Dedekind cut (zero) $(A,B)$ such that for all Dedekind cut $(C,D)$
	\begin{equation}
		(A,B)+(C,D)=(C,D)
	\end{equation}
	
	\item Inverse element by addition exists not for all Dedekind cuts.
	
	\item Sum of proper and improper Dedekind cuts is improper Dedekind cut.
	
	\item Sum of proper Dedekind cuts is proper Dedekind cut.
\end{enumerate}

\textbf{Proof}\\

\begin{enumerate}
	\item Consider two Dedekind cuts $(A,B)$ and $(C,D)$ on $\V{Q}{1}$ such that
	\begin{itemize}
		\item
		\begin{equation}
			S=\{x\in\V{Q}{1}|\forall y\in\V{Q}{0}^+|x|<y\}
		\end{equation}
		
		\item
		\begin{equation}
			\begin{cases}
				A=\{x\in\V{Q}{1}|\forall y\in S\;\;x<y\}\\
				D=\{x\in\V{Q}{1}|\forall y\in S\;\;x>y\}\\
			\end{cases}
		\end{equation}
		
		\item
		\begin{equation}
			\begin{cases}
				B=\V{Q}{1}\setminus (A\cup\sup(A))\\
				C=\V{Q}{1}\setminus (D\cup\sup(D))\\
			\end{cases}
		\end{equation}
	\end{itemize}
	From definition of $A$ and $C$ it follows that $A+C=A$. But from definition of $B$ and $D$ it follows that $B+D=D$. $(A,D)$ is not a Dedekind cut because
	\begin{equation}
		\forall x\in S
		\begin{cases}
			x\notin A\\
			x\notin D\\
		\end{cases}
	\end{equation}
	
	\item As written before, commutativity and associativity follows from corresponding properties of $\V{Q}{\alpha}$.
	
	\item Zero Dedekind cut is $(A,B)$ where
	\begin{equation}
		\begin{cases}
			A=\{x\in \V{Q}{\alpha}|x<0\}\\
			B=\{x\in \V{Q}{\alpha}|x>0\}\\
		\end{cases}
	\end{equation}
	
	\item Consider Dedekind cut $(C,D)$ and set $S$ from first point. Let $(E,F)$ be Dedekind cut such that
	\begin{equation}
		(C,D)+(E,F)=0
	\end{equation}
	If $\sup(E)>0$ then $\sup(C+E)>0$, because $\sup(C)>0$. Thus case where $\sup(E)>0$ does not fit. We can consider three cases:
	\begin{itemize}
		\item $\sup(E)<\inf(S)$, in this case
		\begin{equation}
			\forall x\in E\;\;\forall y\in C\;\;x+y<\inf(S)
		\end{equation}
		hence $C+E$ is not first part of zero Dedekind cut.
		
		\item $\inf(S)<\sup(E)<0$, in this case
		\begin{equation}
			0<\sup(C+E)<\inf(S)
		\end{equation}
		hence $C+E$ is not first part of zero Dedekind cut.
		
		\item $\sup(E)=\inf(S)$, in this case we have situation similar to point 1, i.e. $(C,D)+(E,F)$ is not a Dedekind cut.
	\end{itemize}
	
	\item Let $Q$ be set of Dedekind cuts of $\V{Q}{\alpha}$ and $(A,B)$ be proper Dedekind cut and $(C,D)$ be improper Dedekind cut, now we will check Dedekind cut properties for $(A,B)+(C,D)=(A+C,B+D)$:
	\begin{itemize}
		\item first property is obvious:
		\begin{equation*}
			\forall q\in(A+C)\;\;\forall p\in(B+D)\;\;q<p
		\end{equation*}
		
		\item let us assume that exist two elements $q_1,q_2\in\V{Q}{\alpha}$ such that $\{q_1,q_2\}\subseteq Q\setminus(A+C,B+D)$ and $q_1<q_2$. Thus
		\begin{equation}\label{SomeWRNewProofEq3}
			s\in A+C
			\Rightarrow
			s+(q_2-q_1)\notin B+D
		\end{equation}
		but we know that\footnote{Because $(A,B)$ is proper Dedekind cut and $(C,D)$ is improper Dedekind cut.}
		\begin{itemize}
			\item
			\begin{equation}
				\forall\varepsilon\in\V{Q}{\alpha}^+\;\;\exists\text{ ball }\ball\text{ of radius }\varepsilon:A\cap\ball\ne\varnothing\wedge B\cap\ball\ne\varnothing
			\end{equation}
			
			\item if $E$ is set of all balls such that
			\begin{equation}
				\ball\in E\Rightarrow
				C\cap\ball\ne\varnothing\wedge D\cap\ball\ne\varnothing\
			\end{equation}
			then
			\begin{equation}\label{SomeWRNewProofEq}
				\exists d\in\V{Q}{\alpha}^+:\forall\ball\in E\;\;\sup(\ball)-\inf(\ball)>d
			\end{equation}
		\end{itemize}
		Hence if $M$ is set of all balls such that
		\begin{equation}
			\ball\in M
			\Rightarrow
			(A+C)\cap\ball\ne\varnothing\wedge(B+D)\cap\ball\ne\varnothing
		\end{equation}
		then
		\begin{equation}\label{SomeWRNewProofEq2}
			\forall\ball\in M\;\;\sup(\ball)-\inf(\ball)>d
		\end{equation}
		where $d$ is same as in equation (\ref{SomeWRNewProofEq}).
		
		Thus we have proven that if $(A,B)+(C,D)$ is Dedekind cut, then it is improper.
		
		Let us notice that (it follows from equation (\ref{SomeWRNewProofEq3})) $q_2-q_1>d$ where $d$ is same as in equations (\ref{SomeWRNewProofEq}) and (\ref{SomeWRNewProofEq2}). Hence
		\begin{equation}
			\exists p\in(A+C):p+(q_2-q_1)\in(B+D)
		\end{equation}
		But because
		\begin{equation}
			\exists p\in(A+C):p+d\in(B+D)
		\end{equation}
		then $p+q_2\in(B+D)$, hence we have contradiction and $Q\setminus((A+C)\cup(B+D))$ can be $\varnothing$ or consist of one element.
		
		\item claim that
		\begin{equation}
			\begin{cases}
				\inf(B+D)\notin(B+D)\\
				\sup(A+C)\notin(A+C)\\
			\end{cases}
		\end{equation}
		follows from corresponding properties of $(A,B)$ and $(C,D)$.
		
		\item claim that
		\begin{equation}
			(A+C)\cap(B+D)
		\end{equation}
		follows from first and third properties of Dedekind cut.
	\end{itemize}
	Thus we have proven that sum of proper and improper Dedekind cuts is improper Dedekind cut.
	
	\item Consider proper Dedekind cuts $(A,B),(C,D)\in\V{\overline{R}}{\alpha}$, according to definition of proper Dedekind cuts, for all $\varepsilon\in\V{Q}{\alpha}$ exists balls $\ball_1$ and $\ball_2$ of radius $\frac{\varepsilon}{2}$ such that ball $\ball=\ball_1+\ball_2$ has radius $\varepsilon=\frac{\varepsilon}{2}+\frac{\varepsilon}{2}$. Because $\varepsilon$ can be arbitrary small, then sum of proper Dedekind cuts is proper Dedekind cut.
\end{enumerate}

$\square$

As previously shown, not all Dedekind cuts have inverse element by addition. But let us define $-(A,B)$ for all Dedekind cut $(A,B)$.

\begin{definition}[Opposite Dedekind cut]
	Let $(A,B)\in\V{\overline{R}}{\alpha}$, then Dedekind cut $(E,F)=-(A,B)$ if
	\begin{equation}
		\begin{cases}
			E=\{x\in\V{Q}{\alpha}|\exists y\in B:x=-y\}\\
			F=\{x\in\V{Q}{\alpha}|\exists y\in A:x=-y\}\\
		\end{cases}
	\end{equation}
\end{definition}

\begin{theorem}[About opposite Dedekind cut]\label{OppositeDedCutWRTheor}
	If $x=(A,B)\in\V{\overline{R}}{\alpha}$ is a proper Dedekind cut, then
	\begin{equation}
		x+(-x)=0
	\end{equation}
\end{theorem}

\textbf{Proof}\\

From definition of opposite Dedekind cut it follows that if $(C,D)=x+(-x)$ then
\begin{equation}
	\begin{cases}
		C=\{y\in\V{Q}{\alpha}|y<0\}\\
		D=\{y\in\V{Q}{\alpha}|y>0\}\\
	\end{cases}
\end{equation}
Thus because $(C,D)=0$ then
\begin{equation}
	x+(-x)=0\text{ for all proper }x
\end{equation}

$\square$

\begin{definition}[Absolute value]
	Let $x\in\V{\overline{R}}{\alpha}$, then
	\begin{equation}
		|x|=
		\begin{cases}
			x,\quad x\geqslant 0\\
			-x\quad x<0\\
		\end{cases}
	\end{equation}
\end{definition}

\begin{definition}[Multiplication of positive Dedekind cuts]
	Let $(A,B)>0$ and $(C,D)>0$ be Dedekind cuts on $\V{Q}{\alpha}$, then pair $(E,F)$ is product of $(A,B)$ and $(C,D)$ if
	\begin{equation}
		\begin{cases}
			E=A\cdot C\\
			F=B\cdot D\\
		\end{cases}
	\end{equation}
\end{definition}

\begin{definition}[Multiplication of Dedekind cuts with arbitrary sign]
	Let $(A,B)$ and $(C,D)$ be Dedekind cuts on $\V{Q}{\alpha}$ and
	\begin{equation}
		(E,F)=|(A,B)|\cdot|(C,D)|
	\end{equation}
	then
	\begin{equation}
		(A,B)\cdot(C,D)=
		\begin{cases}
			(E,F),\quad(A,B)\text{ and }(C,D)\text{ have same sign}\\
			-(E,F),\quad\text{otherwise}\\
		\end{cases}
	\end{equation}
\end{definition}

\textbf{Properties of multiplication}\\

\begin{enumerate}
	\item Product of Dedekind cuts is not necessary a Dedekind cut.
	
	\item Product of Dedekind cuts is commutative, associative and distributive with addition (it follows from corresponding properties of multiplication on $\W{Q}{\alpha}$).
	
	\item Exists Dedekind cut (unit) $(A,B)$ such that for all Dedekind cut $(C,D)$
	\begin{equation}
		(A,B)\cdot(C,D)=(C,D)
	\end{equation}
	
	\item Inverse element by multiplication exists not for all Dedekind cuts.
\end{enumerate}

\textbf{Proof}\\

\begin{enumerate}
	\item Consider two Dedekind cuts $(A,B)$ and $(C,D)$ on $\W{Q}{1}$ such that
	\begin{itemize}
		\item
		\begin{equation}
			\begin{cases}
				A=\{x\in\V{Q}{1}|\forall y\in\V{Q}{0}^+\;\;x<y\}\\
				D=\{x\in\V{Q}{1}|\forall y\in\V{Q}{0}^+\;\;x>y\}
			\end{cases}
		\end{equation}
		
		\item
		\begin{equation}
			\begin{cases}
				B=\V{Q}{1}\setminus(A\cup\sup(A))\\
				C=\V{Q}{1}\setminus(D\cup\sup(D))\\
			\end{cases}
		\end{equation}
	\end{itemize}
	Let us notice that
	\begin{equation}
		\begin{cases}
			A\cdot C=A\\
			B\cdot D=D\\
		\end{cases}
	\end{equation}
	Thus
	\begin{equation}
		\forall x\in\V{Q}{0}^+
		\begin{cases}
			x\notin A\cdot C\\
			x\notin B\cdot D\\
		\end{cases}
	\end{equation}
	hence $(A\cdot C,B\cdot D)$ is not Dedekind cut.
	
	\item As written before, commutativity, associativity and distributivity follows from corresponding properties of $\V{Q}{\alpha}$.
	
	\item Unit Dedekind cut is $(A,B)$ where
	\begin{equation}
		\begin{cases}
			A=\{x\in\V{Q}{\alpha}|x<1\}\\
			B=\{x\in\V{Q}{\alpha}|x>1\}\\
		\end{cases}
	\end{equation}
	
	\item Consider Dedekind cut $(A,B)$ from first point. Let $(C,D)$ be Dedekind cut such that
	\begin{equation}\label{MultPropWREq}
		\begin{cases}
			A\cdot C=\{x\in\V{Q}{1}|x<1\}\\
			B\cdot D=\{x\in\V{Q}{1}|x>1\}\\
		\end{cases}
	\end{equation}
	I.e.
	\begin{equation}
		(A,B)\cdot(C,D)=1
	\end{equation}
	
	Because unit and $(A,B)$ are greater then zero, then $(C,D)$ must be greater then zero too.
	
	Let us define set $S$ like in addition case
	\begin{equation}
		S=\{x\in\V{Q}{1}|\forall y\in\V{Q}{0}^+|x|<y\}
	\end{equation}
	
	Because $1>\sup(S)$, then to satisfy the condition in equation (\ref{MultPropWREq}) it is necessary that
	\begin{equation}
		\forall a\in\V{Q}{0}^+\;\;\sup(C)>a
	\end{equation}
	In other situation
	\begin{equation}
		\forall a\in\V{Q}{0}^+\;\;\sup(A\cdot C)<a
	\end{equation}
	But in this situation because $\inf(D)=\sup(C)$ (requirement of Dedekind cut) we have that
	\begin{equation}
		\forall a\in\V{Q}{0}^+\;\;\inf(D\cdot B)>a
	\end{equation}
	Hence $\inf(D\cdot B)>1$ and then $(A,B)\cdot(C,D)\ne 1$.
\end{enumerate}

$\square$

\begin{definition}[Inverse Dedekind cut]
	Let $(A,B)\in\V{\overline{R}}{\alpha}$, then Dedekind cut $(C,D)=\frac{1}{(A,B)}$ if satisfied one of cases
	\begin{enumerate}
		\item if $(A,B)>0$ then
		\begin{equation}
			\begin{cases}
				C=\left\lbrace x\in\V{Q}{\alpha}\left|\exists y\in B:x=\frac{1}{y}\right.\right\rbrace\cup\V{Q}{\alpha}^{\leqslant 0}\\
				D=\V{Q}{\alpha}\setminus(C\cup\sup(C))\\
			\end{cases}
		\end{equation}
		
		\item if $(A,B)<0$ then
		\begin{equation}
			(C,D)=-\frac{1}{|(A,B)|}
		\end{equation}
	\end{enumerate}
\end{definition}

\begin{theorem}[About inverse Dedekind cut]\label{InverseDedCutWRTheor}
	If $x=(A,B)\in\V{\overline{R}}{\alpha}$ is a proper Dedekind cut, then $\frac{1}{x}$ is proper Dedekind cut and
	\begin{equation}
		x\cdot\frac{1}{x}=1
	\end{equation}
\end{theorem}

\textbf{Proof}\\

Let us prove that if for some $S\subset\V{Q}{\alpha}^+$ $\inf(S)$ is proper, then $\sup\left(\frac{1}{S}\right)$ is proper too.

Consider $S\subset\V{Q}{\alpha}^+$ such that $\inf(S)$ is proper. In such case for all $\varepsilon\in\V{Q}{\alpha}^+$ exist $x<\inf(S)$ and $y>\inf(S)$ such that $y-x<\varepsilon$. Also let us consider following set
\begin{equation}
	\frac{1}{S}=D=\left\lbrace x\in\V{Q}{\alpha}\left|\exists y\in S:x=\frac{1}{y}\right.\right\rbrace
\end{equation}
Let us fix $\varepsilon\in\V{Q}{\alpha}^+$ and let $x<\inf(S)$ and $y>\inf(S)$ such that $y-x<\varepsilon$. Consider
\begin{equation}
	\frac{1}{x}-\frac{1}{y}=\frac{y-x}{x\cdot y}
\end{equation}
Because for all $\varepsilon\in\V{Q}{\alpha}^+$ $y-x<\varepsilon$ and $x\cdot y$ is limited from above, then for all $\delta\in\V{Q}{\alpha}^+$ we can choose $x<\inf(S)$ and $y>\inf(S)$ such that
\begin{equation}
	\frac{1}{x}-\frac{1}{y}=\frac{y-x}{x\cdot y}<\delta
\end{equation}
hence $\sup\left(\frac{1}{S}\right)$ is proper too.

Let us assume that $(A,B)>0$ is proper and
\begin{equation}
	C'=\left\lbrace x\in\V{Q}{\alpha}\left|\exists y\in B:x=\frac{1}{y}\right.\right\rbrace
\end{equation}
Using just proven claim we have that $\sup(C')$ is proper too\footnote{$C'=\frac{1}{B}$}. Because $\sup(C)=\sup\left(C'\cup\V{Q}{\alpha}^{\leqslant 0}\right)=\sup(C')$ then $\sup(C)$ is proper.

Let $(A,B)>0$ be proper Dedekind cut. Because $C$ is first set of Dedekind cut $\frac{1}{(A,B)}$ and $\sup(C)$ is proper, then $\frac{1}{(A,B)}=\sup(C)$ is proper.

If $(A,B)<0$ ($(A,B)$ is still proper), then $\frac{1}{|(A,B)|}$ is proper and hence $-\frac{1}{|(A,B)|}=\frac{1}{(A,B)}$ is proper too.

Claim that $(A,B)\cdot\frac{1}{(A,B)}=1$ follows from fact that for all $S_1,S_2\subset\V{Q}{\alpha}^+$ $\inf(S_1)\cdot\inf(S_2)=\inf(S_1\cdot S_2)$.

$\square$

\begin{theorem}[About connection of addition and order on $\V{\overline{R}}{\alpha}$]\label{ConnectAdditionWRTheor}
	Let $a,b,c\in\V{\overline{R}}{\alpha}$, then
	\begin{equation}
		a\geqslant b\Leftrightarrow
		a+c\geqslant b+c
	\end{equation}
\end{theorem}

\textbf{Proof}\\

Let $(A,B),(C,D),(E,F)\in\V{\overline{R}}{\alpha}$ and $(A,B)\geqslant(C,D)$.
According to definitions of order and addition we have
\begin{equation}
	\begin{cases}
		C\subseteq A\\
		(A,B)+(E,F)=(A+E,B+F)\\
		(C,D)+(E,F)=(C+E,D+F)\\
	\end{cases}
	\Rightarrow
	C+E\subseteq A+E
\end{equation}

$\square$

\begin{theorem}[About connection of multiplication and order on $\V{\overline{R}}{\alpha}$]\label{ConnectMultWRTheor}
	Let $a,b\in\V{\overline{R}}{\alpha}$ such that $a\geqslant 0$ and $b\geqslant 0$, then
	\begin{equation}
		a\cdot b\geqslant 0
	\end{equation}
\end{theorem}

\textbf{Proof}\\

Let $(A,B),(C,D)\in\V{\overline{R}}{\alpha}$ and $(A,B)\geqslant 0$, $(C,D)\geqslant 0$. This means that
\begin{equation}
	\begin{cases}
		A\cap\V{Q}{\alpha}^+\ne\varnothing\\
		C\cap\V{Q}{\alpha}^+\ne\varnothing\\
	\end{cases}
\end{equation}
Hence
\begin{equation*}
	(A\cdot C)\cap\V{Q}{\alpha}^+\ne\varnothing
	\Rightarrow
\end{equation*}
\begin{equation*}
	(A,B)\cdot(C,D)\geqslant 0
\end{equation*}

$\square$
	
	\section{$\W{R}{\alpha}$ and $\V{R}{\alpha}$}\label{RealDefSection}
	In this section we will define extensions of real field that are fields too.

\subsection{Definition, order and embedding}

\begin{definition}[$\V{R}{\alpha}$]
	$\V{R}{\alpha}$ is set of all proper elements of $\V{\overline{R}}{\alpha}$ with inherited sum, multiplication and linear order.
\end{definition}

\begin{definition}[$\W{R}{\alpha}$]
	\begin{equation}
		\W{R}{\alpha}=\V{R}{\omega_{\alpha}}
	\end{equation}
\end{definition}

\begin{theorem}[About $\V{R}{\alpha}$]
	$\V{R}{\alpha}$ is linearly ordered field.
\end{theorem}

\textbf{Proof}\\

Commutativity, associativity and distributivity of addition and multiplication on $\V{R}{\alpha}$ are inherited from $\V{\overline{R}}{\alpha}$.

Invertibility of every element of $\V{R}{\alpha}$ relative to addition and multiplication follows from theorems \ref{OppositeDedCutWRTheor} (About opposite Dedekind cut) and \ref{InverseDedCutWRTheor} (About inverse Dedekind cut). From same theorems it follows that $0$ is neutral element of addition and $1$ is neutral element of multiplication.

Connection of binary operations with order inherits from theorems \ref{ConnectAdditionWRTheor} and \ref{ConnectMultWRTheor}.

Hence, $\V{R}{\alpha}$ is linearly ordered field.

$\square$

Embedding $\V{Q}{\alpha}\subset\V{R}{\alpha}$ is given by definition \ref{QREmbedDef}. (Let us notice that every rational element is proper).

\subsection{Topology}

Topology on $\V{R}{\alpha}$ is inherited from $\V{\overline{R}}{\alpha}$.

\begin{theorem}[About disconnection of $\V{R}{\alpha}$]
	If $\alpha>0$ then $\V{R}{\alpha}$ is disconnected in topology induced from $\V{\overline{R}}{\alpha}$.
\end{theorem}

\textbf{Proof}\\

Consider set
\begin{equation}
	S=\{x\in\V{R}{\alpha}|\forall y\in\mathbb{R}\;\;|x|<|y|\}
\end{equation}

Set $S$ has empty border in $\V{R}{\alpha}$ because its border consists of improper elements of $\V{\overline{R}}{\alpha}$, hence it does not belong to $\V{R}{\alpha}$. Thus $S$ is opened and closed simultaneously, hence $\V{R}{\alpha}$ is disconnected.

$\square$

\begin{definition}[Ball (opened ball) on $\V{R}{\alpha}$]
	Ball on $\V{R}{\alpha}$ with center in $x_0\in\V{R}{\alpha}$ and radius $r\in[0,+\infty)_{\alpha}$ is set $\ball_{\alpha}(x_0,r)$ such that
	\begin{equation}
		\ball_{\alpha}(x_0,r)=\{x\in\V{R}{\alpha}|r>|x-x_0|\}
	\end{equation}
\end{definition}

Let us notice that $\ball_{\alpha}(x_0,0)=\varnothing$.

Because for every ball $B\subset\V{\overline{R}}{\alpha}$ its intersection $B\cap\V{R}{\alpha}$ can be represented through union of opened balls on $\V{R}{\alpha}$, then set of all opened balls on $\V{R}{\alpha}$ is base of topology inducted from $\V{\overline{R}}{\alpha}$.
	
	\printbibliography
\end{document}